\documentclass[12pt]{amsart}
\usepackage{amsmath}

\def\e{\epsilon}
\def\lf{\left}
\def\ri{\right}

\def\wt{\widetilde}

\def\p{\partial}

\newcommand\ce{{\mathbb C}}
\newcommand\C{{\mathbb C}}

\def\ii{\sqrt{-1}}
\def\jbar{{\bar\jmath}}

\def\ttg{\tilde{g}}

\def\tg{\tilde{g}_{i\bar\jmath}}

\def\G{g_{ij}}

\def\K{K\"ahler }
\def\KR{K\"ahler-Ricci }
\def\KRF{K\"ahler-Ricci flow }
\def\KRS{K\"ahler-Ricci soliton }
\def\KRS{K\"ahler-Ricci soliton }

\def\tM{\wt M}
\def\tN{\wt N}
\def\tL{\wt L}

\def\be{\begin{equation}}
\def\ee{\end{equation}}

\def\lf{\left}
\def\ri{\right}

\def\e{\epsilon}

\def\ijb{{i\jbar}}

\def\Ric{\text{Ric}}

\def\wt{\widetilde}

\def\p{\partial}

\def\C{\Bbb C}
\def\cn{\Bbb C^n}

\def\wt{\widetilde}

\def\p{\partial}

\def\C{\Bbb C}
\def\ii{\sqrt{-1}}
\def\ttg{\wt{g}}

\newtheorem{thm}{Theorem}[section]

\newtheorem{lem}{Lemma}[section]

\newtheorem{cor}{Corollary}[section]
\theoremstyle{definition}

\theoremstyle{remark}
\newtheorem{rem}{Remark}
\numberwithin{equation}{section}
 \begin{document}
\author{Albert Chau$^1$}
\thanks{$^1$Research
partially supported by NSERC grant no. \# 327637-06}

\address{Waterloo University, Department of Pure Mathematics,
  200 University avenue, Waterloo, ON N2L 3G1, CANADA}
\email{a3chau@math.uwaterloo.ca}

\author{Luen-Fai Tam$^2$}
\date{Febuary,  2007}
\thanks{$^2$Research
partially supported by Earmarked Grant of Hong Kong \#CUHK403005}
\address{Department of Mathematics, The Chinese University of Hong Kong,
Shatin, Hong Kong, China.} \email{lftam@math.cuhk.edu.hk}

\title{A survey on the K\"ahler-Ricci flow and Yau's uniformization conjecture}

\begin{abstract}
 Yau's uniformization conjecture states: a complete
 noncompact \K  manifold  with positive holomorphic bisectional curvature
 is biholomorphic to $\ce^n$.  The \KRF has provided a powerful
 tool in understanding the conjecture, and has been used to verify
 the conjecture in several important cases.
 In this article we present a survey of the \KRF with focus on its application to
 uniformization. Other interesting methods and results related to
 the study of    Yau's conjecture are also discussed.
\end{abstract}

\maketitle \markboth{Albert Chau and Luen-Fai Tam} {A survey on the K\"ahler-Ricci flow and Yau's uniformization conjecture}

\section{Introduction}

A fundamental problem in complex geometry is to generalize the
classical uniformization theorems on Riemman surfaces to higher
dimensions.  In \K geometry, the problem is to determine how
curvature affects the underlying holomorphic structure of a \K
manifold.  In one complex dimension it is well known
that a complete simply connected Riemann surface $(M, g)$ is biholomorphic to either the Riemann sphere (when $M$ is compact) or the
complex plane (when $M$ is noncompact) if the curvature is
positive, and it is biholomorphic to the open unit disc if the
curvature is negative and bounded from above away from zero.

In higher dimensions, for the compact case, the famous conjecture
of
 Frankel says that a compact K\"ahler manifold with
positive holomorphic bisectional curvature is biholomorphic to
$\mathbb{CP}^n$.  Frankel's conjecture was proved by Siu-Yau
\cite{SY}.  The stronger Hartshorne conjecture was proved by Mori
\cite{M}.   In case of    compact K\"ahler manifolds with
nonnegative holomorphic bisectional curvature the uniformization
of such manifolds was determined by Bando \cite{Bando84} for
complex dimension three and by Mok \cite{Mok1988} for all
dimensions.

In the complete noncompact case there is a long
standing conjecture due to Yau \cite{Y} in 1974 predicting analogous results
 for positively curved non-compact \K manifolds:

{\bf Yau's Conjecture} \cite{Y} {\it A complete noncompact
\K manifold with positive holomorphic bisectional curvature
is biholomorphic to $\cn$. }

Or in Yau's words \cite[p.620]{Y}: {\it The question is to
demonstrate that every   noncompact K\"ahler manifold with
positive bisectional curvature is biholomorphic to the complex
eulcidean space.}

 The conjecture continues to generate much research activity,
and although the full conjecture remains unproved, this research
has produced many nice results and useful techniques and methods.
In this survey we will discuss the method of using evolution
equations, in particular the \KR flow, to study the conjecture.

 In \cite{Sh}, Shi (under the supervision of Yau) began a program to use the \KRF equation
\begin{equation}\label{s0e1}
   \frac{d g_{i\jbar}}{dt}=-R_{i\jbar}
\end{equation}
to prove Yau's uniformization conjecture \footnote{The \KR flow had been used in Bando's work \cite{Bando84}
(under the supervision of Yau) and Mok's work \cite{Mok1988} on
the uniformization of compact \K manifolds with
nonnegative holomorphic bisectional curvature.  The use of the \KRF in these works however was not very extensive.}.  (\ref{s0e1}) is a
geometric  evolution equation which deforms an initial \K metric
in the direction of its Ricci curvature (see \S2).  The equation
is a strictly parabolic system for $g_{i\jbar}$, so one generally
expects the geometry to improve under this evolution.  The idea is
to show that the geometry improves to the point that it determines
the holomorphic structure of the manifold.  This idea has been
successfully applied in many cases now.  In this article we survey
the progress on Yau's conjecture and the \KRF program.  We will
also survey some more general results relating to the conjecture.
Since there are many works, the survey cannot be exhaustive.

 An outline of the paper is as follows.  In \S 2 we introduce
 the \KRF and its basic theory on complete noncompact \K manifolds,
 and on nonnegatively curved \K manifolds in particular.
 In \S 3 we discuss the rate of curvature decay on a complete noncompact
 nonnegatively curved \K manifold, which is  also related to  volume growth. This turns out to be useful for latter application.
  The Steinness of non-negatively curved noncompact
  \K manifolds is discussed in \S 4, and this begins our actual survey
  of uniformization results.  In \S 5 we introduce the \KR solitons
  and present our uniformization theorem for steady
  and expanding gradient \KR solitons.
  We place this section here because \KR solitons are canonical
   models for the cases in \S 6 and the analogy drawn there between
   the \KRF and complex dynamical systems.
    In \S 6 the connection between the \KRF and complex dynamical systems
    is developed  and used to prove a uniformization theorem
      for eternal solutions to the \KRF and normalized \KR flow.
      As a corollary, we state a uniformization theorem
        for nonnegatively curved \K manifolds of average quadratic curvature
        decay, in particular for manifolds with maximum volume growth.
        Finally, in \S 7 we present the gap theorems for nonnegatively
        curved \K manifolds of faster
        than quadratic curvature decay. Whenever it is possible, we will
        sketch the main ideas of the proofs of the results.

 \section{\KRF}
The Ricci flow on a complete Riemannian manifold $(M, g)$ is
the following evolution equation for the metric $g$:
\begin{equation}\label{s1e0}
   \left\{%
\begin{array}{ll}
    \frac{d g(t)}{dt}&=-2Rc(t) \hbox{ } \\
    g(0)&=g\hbox{,} \\
\end{array}%
\right.
\end{equation}
where $Rc(t)$ is the Ricci tensor of $g(t)$. The Ricci flow was
introduced by Hamilton \cite{Ha} to study the Poincar\'e
conjecture. In this seminal work, Hamilton proved (\ref{s1e0}) has
a short time solution on any smooth compact Riemannian manifold.
Hamilton then showed that the solution exists for all time and
converges after rescaling on any compact three manifold with positive Ricci
curvature.

 We are interested in complete noncompact Riemannian manifolds. The first fundamental result for the Ricci flow on noncompact manifolds is the following short time existence result of Shi \cite{Sh0}:
\begin{thm}\label{Shi-shorttime}\cite{Sh0} Let $(M^n,g)$ be a complete
noncompact Riemannian manifold with bounded sectional curvature.
Then there exists $0<T<\infty$, depending only on the initial
curvature bound, such that (\ref{s1e0}) has a solution $g(t)$ on
$M\times [0, T]$.  Moreover, for all $t\in [0, T]$ we have
\begin{enumerate}
\item[(i)] $g(t)$ has bounded sectional curvature and is
equivalent to $g$.
 \item[(ii)] For any integer $m\ge0$, there is a constant
 depending only on $n, m$ and the bound of the   curvature of the
 initial metric
 $g$ such that
 $$
 \sup_{x\in M}|\nabla^mRm|(x,t)\le \frac{C}{t^m}
 $$
for all $0\le t\le T$. Here $\nabla$ is the covariant derivative
with respect to $g(t)$.
\end{enumerate}
\end{thm}

The following theorem of Shi \cite{Sh2} is fundamental to the
theory and application of the Ricci flow to \K manifolds.

 \begin{thm}\label{Shi-shorttimeKahler}\cite{Sh2}  Let $(M,g)$ be as in Theorem
 \ref{Shi-shorttime}. Suppose $(M,g)$ is K\"ahler.  Then
\begin{enumerate}
\item [(i)] the solution $g(t)$ in the Theorem \ref{Shi-shorttime}
  is  \K   for $t\in [0,T]$; and
\item[(ii)] $g(t)$ has non-negative holomorphic bisectional
curvature  if the same is true for $g$.
\end{enumerate}
\end{thm}
Hence in the K\"ahler category, we refer to the Ricci flow as the
\KR flow and the equation becomes:
\begin{equation}\label{s1e1}
   \frac{d g_{i\jbar}}{dt}=-R_{i\jbar}.
\end{equation}

For compact K\"ahler manifolds, part (i) of Theorem
\ref{Shi-shorttimeKahler} was proved by Bando \cite{Bando84} and
part (ii) was proved by Bando \cite{Bando84} for complex three
manifolds then by Mok \cite{Mok1988} for all dimensions. Shi's
proof of (i) for the noncompact case is different from that in
\cite{Bando84}. Shi's proof of (ii) is similar to that in
\cite{Bando84} and \cite{Mok1988} which uses a so-called
null-vector condition introduced in \cite{Bando84}.  The proofs of
both (i) and (ii) use  the maximum principle for noncompact
Riemannian manifolds which relies on the following result which
was proved by Shi \cite{Sh2} using  an idea of Greene-Wu
\cite{GW76}: Let $(M,g)$ be a complete noncompact manifold with
bounded curvature. Then one can construct a smooth function with
bounded gradient and bounded Hessian, which is uniformly
equivalent to the distance function from a point.

 We now focus on complete non-compact \K manifolds with non-negative
  holomorphic bisectional curvature.  The main long time existence result
  for the \KRF in this setting is the following theorem
   of Shi \cite{Sh2}.  We will denote the scalar curvature by $\mathcal{R}$ and define the average scalar curvature $k(x, r)$ as:
$$k(x, r):=\frac{1}{V_x(r)}\int_{B_x(r)}\mathcal{R} dV.$$

\begin{thm}\label{Shi-longtime}\cite{Sh2} Let $(M^n, g)$ be a
complete noncompact K\"ahler manifold with bounded and nonnegative
holomorphic bisectional curvature. Suppose
\begin{equation}\label{longtimedecay}
 k(x,r)\le \frac{C}{
(1+r)^\e}
\end{equation} for some constants $C$ and $0< \e\le 2$, and all $x$
and $r$. Then (\ref{s1e1}) has a long time solution $g(t)$.
Moreover,
\begin{enumerate}
\item [(i)] $$ \mathcal{R}(x,t)\le Ct^{\frac{ 2(1-\e)}{\e}},$$ for
some constant $C$ for all $x$ and $t$.

\item[(ii)] For any integer $m\ge0$, there is a constant $C$ such
that
$$|\nabla^mRm|(x,t)\le Ct^{\frac{ 2(1-\e)}{\e}(m+2)} $$
for all $x$ and $t$.
\end{enumerate}
\end{thm}
For example, if $\e=2$, then $t\mathcal{R}$ will be uniformly
bounded on spacetime.

 The proof of the theorem is based on estimating the
volume element $$F(x,t)=\log \frac{\det(g_{\ijb}(x,t))}{
\det(g_{\ijb}(x,0))}.$$   Shi proved that $F$ stays bounded from
below  on any finite time interval.  Note that the bound $F\le0$
from above follows from Theorem \ref{Shi-shorttimeKahler}.  Using
this and a parabolic version of the third order estimate for the
complex Monge-Amp\`ere equation to conclude that \KR flow, Shi
then showed that the \KRF cannot develop a singularity in finite
time.    In the case of $\e>1$, there is also a method by Ni-Tam
\cite{NT} for proving the long time existence.  Namely, by
extending the method of Mok-Siu-Yau \cite{MSY}, Ni-Shi-Tam
\cite{NST} were able to construct a potential $u_0$, for the Ricci
form of the initial metric, having uniformly bounded gradient.  It
is readily seen that $u=u_0-F$ is a potential of the Ricci form
for $g(t)$.  Using the maximum principle, one can show that $u$
also has uniformly bounded gradient.  One then shows that $|\nabla
u|^2+ \mathcal{R}$ is a subsolution of the time dependent heat
equation, and applying the maximum principle again, one concludes
that $|\nabla u|^2+ \mathcal{R}$ is bounded by its maximum at
$t=0$.  Since $g(t)$ has nonnegative holomorphic bisectional
curvature, the Riemannian curvature can also be bounded from this.

On the other hand, the estimate of $F$ depends on the
fact that
$$
\Delta_0 F\le \mathcal{R}(0)+e^F\frac{\p F}{\p t}
$$
and some mean value inequalities (see \cite{NT} for example). Here
$\Delta_0$ is the Laplacian of the initial metric. Once we estimate
$F$ we can then estimate $ \mathcal{R}$ using the fact that
$$
F(x,t)=-\int_0^t\mathcal{R}(x,\tau)d\tau
$$
and the Li-Yau-Hamilton type Harnack inequality Theorem
\ref{caoharnack} of Cao \cite{cao} which implies that
$t\mathcal{R}$ is nondecreasing in time.  In any case, the estimates of the
covariant derivatives of the curvature tensor follow from the
general method developed by Shi \cite{Sh0}.

Note that in the theorem, condition (\ref{longtimedecay}) is
uniform in $x$.  In many cases this condition is at least true at a
point by \cite{NT2}. Using this and the pseudolocality of Ricci
flow by Perelman \cite{P1}, which can be generalized to complete
noncompact manifolds, Yu and the authors \cite{ChauTamYu2007} obtained the following result on long time existence:

\begin{thm}\label{ChauTamYu}
Let $(M^n, g)$ be a complete non-compact \K  \ manifold with
non-negative holomorphic bisectional curvature with injectivity
radius bounded away from zero such that $$|Rm|(x)\to0$$ as $x\to
\infty$.     Then   the \KR flow with initial data $g$ has a long
time solution $g(t)$ on $M\times [0, \infty)$.
\end{thm}
It would be nice if one can remove the condition on the
injectivity radius.

 \section{Curvature decay rate}
 In order to apply   \KRF to Yau's conjecture, the following question becomes important in light of Theorem \ref{Shi-longtime}:
 What can we expect of the the curvature decay rate on complete
 noncompact K\"ahler manifolds with positive bisectional curvature?
 There are some results on the volume and curvature of manifolds with
 nonnegative holomorphic bisectional curvature in \cite{NST}. For example,
 it was proved that if the scalar curvature of a complete noncompact \K
 manifold with nonnegative holomorphic bisectional curvature decays quadratically
 in the pointwise and average sense, and if the Ricci curvature is positive
 at some point, then the manifold must have maximum volume growth.
 On the other hand, the following was proved by
  Chen-Zhu \cite{Chen-Zhu05}:

 \begin{thm}\label{ChenZhu-curvature}\cite{Chen-Zhu05}  Let $(M^n,g)$ be a complete noncompact K\"ahler manifold with
positive holomorphic bisectional curvature. Then for any $x\in M$
there is a constant
 $C$ which may depend on $x$ such that
 \begin{equation}\label{lineardecay1}
   \frac{1}{V_x(r)}\int_{B_x(r)}\mathcal{R}dV\le \frac{C}{1+r}.
\end{equation}
\end{thm}

Later Ni and Tam \cite{NT2} obtained the following, which
generalizes the above results:
\begin{thm}\label{NiTam-curvature}\cite{NT2}
Let $(M^n,g)$ be a complete noncompact K\"ahler manifold with
nonnegative holomorphic bisectional curvature.
 \begin{itemize}
    \item [(i)] Suppose $M$ is simply connected, then
$M=N\times M'$ holomorphically and isometrically, where $N$ is a
compact simply connected K\"ahler manifold, $M'$ is a complete
noncompact K\"ahler manifold and both $N$ and $M'$ have
nonnegative holomorphic bisectional curvature. Moreover, the
scalar curvature $\mathcal{R'}$ of $M'$ satisfies the linear decay
condition (\ref{lineardecay1}).
    \item [(ii)] If the holomorphic bisectional  curvature of $M$ is
positive at some point, then $M$ itself satisfies the linear decay
condition (\ref{lineardecay1}).
\end{itemize}
\end{thm}
The idea of the proof is to construct a strictly plurisubharmonic
function $u$ of linear growth. If such a $u$ exists, then one can
use $L^2$ theory to construct nontrivial holomorphic section $s$
of the canonical line bundle using $u$ as a weight function. The
growth rate of $u$ will give an estimate of the growth rate of the
$L^2$ norm of the length $||s||$ of  $s$. Using the Bochner type
differential inequality:
\begin{equation}\label{lineardecay2}
\Delta \log ||s||^2\ge \mathcal{R}
\end{equation}
at the points where $s\neq0$, one can show that $\log
||s||$ is at most linear growth by the mean value inequality of Li
and Schoen \cite{LiSchoen84}. Then the growth rate of $\log ||s||$
will give an estimate of the average of $\mathcal{R}$ over
geodesic balls because of (\ref{lineardecay2}).

If the holomorphic bisectional curvature is positive, one just considers the Busemann function which is strictly
plurisubharmonic at a point by a well known result of Wu \cite{W1}. In
general, one may solve the heat equation with the Buesmann
function as initial data.  Then by a careful study of the solution
 for $t>0$, one obtains a suitable strictly plurisubharmonic
 function in the case of Theorem \ref{NiTam-curvature} (ii).

Hence in some sense, linear decay is the slowest decay rate for
the curvature of complete noncompact K\"ahler manifolds with
nonnegative holomorphic bisectional curvature.

In case that the manifold has maximum volume growth, one may
expect the curvature will decay faster. In fact, it was
conjectured by Yau \cite{Y}: If $M$ has maximal volume
 growth in the sense that $V_p(r)\ge Cr^{2n}$ for some $C>0$ for
 some $p\in M$ for all $r$, then the curvature must decay
 quadratically in the average sense. Assuming the curvature is
 bounded, this was confirmed by Chen-Tang-Zhu \cite{CTZ} for
 complex surfaces and Chen-Zhu \cite{Chen-Zhu05} for higher
 dimension under a much stronger assumption that the manifold has
 nonnegative curvature operator. Finally, Ni \cite{Ni05} proved
 this conjecture of Yau in general. More precisely, Ni
 obtained the following:
\begin{thm}\label{Ni-curvature}\cite{Ni05} Let  $(M^n,g)$ be a complete
noncompact K\"ahler manifold with bounded nonnegative
holomorphic bisectional curvature. Suppose $M$ has maximum
volume growth. Then the scalar curvature $\mathcal{R}$
decay quadratically. That is to say, there exists $C>0$
such   that
\begin{equation}\label{quadraticdecay1}
    \frac{1}{V_x(r)}\int_{B_x(r)}\mathcal{R}\,dV\le
    \frac{C}{(1+r)^2}
\end{equation}
for all $x\in M$ and for all $r>0$.
\end{thm}
The main steps of the proof are as follow: First, use the results
in \cite{NT2}, in particular Theorem \ref{NiTam-curvature}, to
prove that a non-flat gradient shrinking K\"ahler-Ricci soliton
with nonnegative bisectional curvature must have zero asymptotic
volume ratio. A gradient shrinking K\"ahler-Ricci soliton is a
K\"ahler metric $g_\ijb$ satisfying
$$
R_\ijb-\rho g_\ijb=f_\ijb, \hspace{12pt} f_{ij}=0
$$
for some smooth real-valued function $f$ and for some $\rho>0$.
For an $n$ dimensional Riemannian manifold with nonnegative Ricci
curvature, the asymptotic volume ratio is defined as:
$$
\mathcal{V}(g) =\lim_{r\to\infty}\frac{V_x(r)}{r^n}.
$$
The limit exists and is independent of the base point $x$
by the Bishop volume comparison theorem.  In the case of
maximal volume growth, this limit is non-zero.

  The second step is
    to prove that any non-flat complete
ancient solution of the \KRF on a K\"ahler manifold with bounded
and nonnegative holomorphic bisectional curvature must also have
zero asymptotic volume ratio for all $t$. This is accomplished by
showing that if this is not true, then by a   blow down argument
as in Perelman \cite{P1} for Riemannian manifolds with nonnegative
curvature operator, one can construct a non-flat gradient
shrinking K\"ahler-Ricci soliton with nonnegative bisectional
curvature which has {\it nonzero}  asymptotic volume ratio.

 Now suppose $(M,g)$
satisfies the conditions in Theorem \ref{Ni-curvature}, then one
can solve the \KRF equation with solution $g(t)$ with initial
condition $g$, and can prove that
$\mathcal{V}(g(t))=\mathcal{V}(g)$ for all $t$. One then uses the second step above to prove that $g(t)$ must exist for all time and
the scalar curvature must satisfy
$$
\mathcal{R}(x,t)\le \frac{C}{1+t}
$$
for some $C$ for all $x$ and $t$.  This is proved by
contradiction: if this were not so, one could construct an ancient
solution as above having non-zero asymptotic volume ratio, which
contradicts the previous assertion.  Finally, one may use this
asymptotic behavior of $\mathcal{R}$ as $t\to\infty$ to get
information of $\mathcal{R}(x,0)$ as $x\to\infty$ and obtain
(\ref{quadraticdecay1}). Here an argument similar to that in
Perelman \cite{P1} is also used.

In order to get more information on the question on curvature
decay rate, it would be helpful to construct examples.  This is considerably easier in the Riemannian setting: constructing complete noncompact Riemannian manidfolds having positive sectional curvature.  It is not easy to construct complete K\"ahler metrics on $\cn$ with positive holomorphic bisectional curvature. The first example
is by Klembeck \cite{Klembeck77} which has positive holomorphic bisectional
curvature.  Klembeck's example has linear decay curvature
with volume growth like $V(r)\sim r^n$.  Later, Cao \cite{Cao94,Cao} constructed examples of $U(n)$ invariant
metrics on $\cn$ having positive holomorphic bisectional curvature.
The examples in \cite{Cao94} have linear curvature decay and volume growth like
$V(r)\sim r^n$ while those in \cite{Cao} have  quadratic curvature decay and maximum volume
growth.

  Complete $U(n)$ invariant K\"ahler metrics with positive holomorphic
  bisectional curvature on $\cn$ have been classified and examples have been constructed by Wu
  and Zheng \cite{WuZheng-pp}.  They have developed a systematic
  way to construct examples having positive holomorphic bisectional
  curvature such that the scalar curvature $\mathcal{R}$ satisfies
  various decay rates.  In particular, for any $1\le\theta\le 2$ there are examples such that

  \begin{equation}\label{WuZheng}
    \frac{1}{V_x(r)}\int_{B_x(r)}\mathcal{R}\,dV\le
    \frac{C}{1+r^\theta}
\end{equation}
for some $C$ for all $r$ and for all $x$.

Finally, there are interesting results by Chen-Zhu
\cite{Chen-Zhu05} on volume growth of positively curved noncompact \K manifolds. For example, they proved that the volume
  of geodesic balls $B_p(r)$  in a complete noncompact
\K manifold with positive holomorphic bisectional curvature
must grow at least like $r^n$, where $n$ is the complex
dimension of the manifold. One may compare this with the
well-known result of Calabi and Yau (see \cite{Y1}, for
example) that in case of Riemannian manifold with
nonnegative Ricci curvature, the growth rate of geodesic
balls is at least linear.

 \section{Steinness of nonnegatively curved manifolds}
 In \cite{Siu}, Siu  asked whether a complete noncompact \K manifold
 with positive holomorphic biectional is Stein.
 This is a very interesting problem, and an affirmative answer to this
  should strongly support Yau's conjecture. In fact, Greene and Wu \cite{GW76} first proved the following result on the
complex structure of complete noncompact K\"ahler manifolds with
positive curvature, which is part of the motivation behind Yau's
conjecture:
\begin{thm}\label{Greene-Wu1} \cite{GW76} Let $(M,g)$ be a complete noncompact
K\"ahler manifold with positive sectional curvature. Then
$M$ is Stein.
\end{thm}
The idea is to produce a smooth strictly plurisubharmonic
exhaustion function.  Then by a well-known result of Grauert, $M$
will be Stein. As in the study of Riemannian manifolds with
positive or nonnegative sectional curvature by Cheeger and Gromoll
\cite{CG72}, the method of Greene-Wu uses the Busemann function
$\mathcal{B}(x)$.  Here the Busemann  function is an exhaustion
function because the sectional curvature is nonnegative
\cite{CG72}. It is also plurisubharmonic if the holomorphic
bisectional curvature is nonnegative and strictly plurisubharmonic
at the points where the holomorphic bisectional curvature is
positive by   Wu \cite{W1}. However, $\mathcal{B}(x)$ is only
Lipschitz, and so one needs to approximate $\mathcal{B}(x)$ by a
smooth strictly plurisubharmonic exhaustion function. This can be
done using the fact that $\mathcal{B}(x)$ is strictly
plurisubharmonic and   the method developed by Greene-Wu
\cite{GW76}. Because of these, Theorem \ref{Greene-Wu1} can be
improved. For example, the theorem is still true if we only assume
that $M$ has nonnegative sectional curvature and positive
holomorphic bisectional curvature.

 In the case that $\mathcal{B}$ is only plurisubharmonic,
 it is more efficient to use another method to
approximate $\mathcal{B}(x)$. Namely, we solve the heat equation
with initial data $\mathcal{B}(x)$.  Ni and Tam \cite{NT2} proved
that if $u(x,t)$ is the solution to this, then $u(x,t)$ is still
plurisubharmonic for $t>0$. Moreover, they showed that the kernel
$\mathcal{K}(x,t)$ of the complex Hessian of $u(x,t)$ is a
parallel distribution. Using this one can obtain the following:

\begin{thm}\label{NiTam-Stein}\cite{NT2}
Let $(M,g)$ be a complete noncompact K\"ahler manifold with
nonnegative holomorphic bisectional curvature. Suppose either $M$
has maximum volume growth or $M$ has a pole.  Then $M$ is Stein.
\end{thm}
Here one uses the fact that $\mathcal{B}(x)$ is still an
exhaustion function by \cite{Shen1996} and so $u(x,t)$ is
also an exhaustion function. As in the proof of Theorem
\ref{NiTam-curvature}, one can   conclude that $u$ is also
strictly plurisubharmonic for $t>0$. In the proof, one
needs the following fact \cite{NT2}: If $(M,g)$ is a
complete K\"ahler manifold with nonnegative holomorphic
bisectional curvature which supports a nontrivial linear
growth harmonic function, then the universal cover of $M$
is a holomorphically isometric to the product of the
complex $\mathbb{C}$ and another complete K\"ahler manifold
with nonnegative holomorphic bisectional curvature.

Without the assumption that the manifold has maximum volume
growth, if the Busemann function is an exhaustion function, then
one can conclude that the universal cover is a production of a
compact Hermitian symmetric manifold and a Stein manifold
\cite{NT2}. Using this, Fangyang Zheng (see \cite{NT2}) obtained
the following:
\begin{thm}\label{Zheng}
Let $(M,g)$ be a complete noncompact K\"ahler manifold with
nonnegative sectional curvature. Then its universal cover
is of the form $\tM=\C^k\times \tN\times \tL $ where $\tN$
is a compact Hermitian symmetric manifold, $\tL$ is Stein
and  $\tL$ contains no Euclidean factor. Moreover,  $M$ is
a holomorphic and Riemannian fiber bundle with fiber
$\tN\times \tL$ over a flat K\"ahler manifold
${\C}^k/\Gamma $. If in addition, the Ricci curvature is
positive at some point, then $M$ is simply connected and
$M=\tM=\tN\times\tL$ where $\tL$ is diffeomorphic to the
Euclidean space.
\end{thm}

 So far the above results have made no assumptions on the boundedness
 of curvature.  When the curvature is bounded, one can use
 the more powerful   \KRF (\ref{s1e1}), which  can be solved for a short
 time by Theorem \ref{Shi-shorttime}. The results in \S2
  suggests the study of $g(t)$ in connection to Yau's conjecture or
  more generally, the Steinness of $M$.
   When $g(t)$ exists for all time $0\leq t<\infty$ we expect the
    asymptotics of $g(t)$ to be particularly useful.
    In \cite{CT5} the authors proved the following:

\begin{thm}\label{fastlineardecay1}\cite{CT5}
 Let $(M^n, g_0)$ be a complete non-compact \K manifold with
bounded non-negative holomorphic bisectional curvature.  Suppose
the scalar curvature of $g_0$ is such that  $k(x, r) \leq k(r)$
for some function $k(r)$ satisfying
\begin{equation}\label{lineardecay}
k(r) \leq \frac{C}{r}
 \end{equation}
as $r\to\infty$ for some $C>0$.  Then $M$ is holomorphically covered by a
pseudoconvex domain in $\C^n$ which is homeomorphic to
$\mathbb{R}^{2n}$.  Moreover, if $M$ has positive
bisectional curvature and is simply connected at infinity, then M
is biholomorphic to a pseudoconvex domain in $\cn$ which is
homeomorphic to $\mathbb{R}^{2n}$, and in particular, $M$ is Stein.
 \end{thm}

When $k(r)=\frac{C}{1+r^{1+\epsilon}}$, for $\epsilon>0$, the
result that $M$ is biholomorphic to a pseudoconvex domain was
proved by Shi \cite{Sh2} under the additional assumption that $(M,
g)$ has positive sectional curvature. Under the same decay
condition and assuming maximum volume growth, similar results were
obtained by Chen-Zhu \cite{Chen}. The  condition of positive
sectional curvature  in \cite{Sh2} was used to produce a convex
compact exhaustion of $M$.  The maximum volume growth condition in
\cite{Chen} was used to control the injectivity radius for $g(t)$.
Note that the decay condition in Theorem \ref{fastlineardecay1} is
almost optimal because of Theorem \ref{NiTam-curvature}.

The idea of proof of Theorem  \ref{fastlineardecay1} is as
follows. Let $g_0$ be as in the theorem. Then one can find $g(t)$
which solves the K\"ahler-Ricci flow equation (\ref{s1e1}) with
initial data $g(0)=g_0$ satisfying the conclusion in Theorem
\ref{Shi-longtime}.  By taking the universal cover and by using
the result of Cao \cite{cao1}, we may assume that $g(t)$ has
positive Ricci curvature.  Let $g_i=g(i)$. Then for any fixed $p\in M$,

\begin{enumerate}
     \item  [\textbf{ (a1)} ] $cg_{i}\le g_{i+1}\le g_i$  for some
     $1>c>0$ for all $i$.

     \item  [\textbf{ (a2)} ] $|\nabla Rm(g_i)|+|Rm(g_i)|\le c'$ for
     some $c'$ on $ B_i(p,r_0)$,  and for some  $r_0>0$  for all $i$
         where  $B_i(p,r_0)$ is the geodesic ball around $p$ with respect to
$g_i$.

 \item  [\textbf{ (a3)} ]  $g_i$ is contracting in the
following sense:   For any $\e$,  for any $i$, there
$i'>i$ with
$$
g_{i'}\le \e g_i
$$
in $B_i(p, r_0)$.
\end{enumerate}
\textbf{(a1)-(a2)} follow from Theorem \ref{Shi-longtime}, and
\textbf{(a3)} follows from Theorem \ref{TH} in \S6.

 Hence Theorem \ref{fastlineardecay1}
is a consequence of the following:
\begin{thm}\label{ChauTam-Stein2}\cite{CT5}  Let $M^n$ be a complex
noncompact manifold. Suppose there exist a
 sequence of complete K\"ahler metrics $g_i$, $i\ge1$ on $M$ with
 properties \textbf{(a1)--(a3)}.
  Then  $M$ is covered by a pseudoconvex
domain in $\cn$ which is homeomorphic to $\mathbb{R}^{2n}$.
\end{thm}
We now sketch some key ideas in the proof Theorem
\ref{ChauTam-Stein2} .  Fix some point $p\in M$. By \textbf{(a2)}
we may lift the metric $g(t)$ to the tangent space via the
exponential map at $p$ and use $L^2$ theory to construct a local
biholomorphism $\Phi_i: D(r)\to M$ for each $i$ so that
$\Phi_i^*(g_i)$ is equivalent to the Euclidean metric in $D(r)$
and equal to the Euclidean metric at $0$, where $D(r)$ is the
Euclidean ball of radius $r$ with center at the origin  (see
Proposition 2.1 in \cite{CT4}).  Thus $\Phi_i$ provides a
holomorphic normal  coordinate at $p$, and one would like to
consider an  appropriate change of coordinate map from  $\Phi_i$
and $\Phi_{i+1}$.  It is not clear one can do this however as
$\Phi_i$ is generally not injective and may  not even be a
covering. Nevertheless, using \textbf{(a1)} and \textbf{(a2)}, we
can find $F_{i+1}:D(r)\to \cn$ such that $\Phi_i=\Phi_{i+1}\circ
F_{i+1}$ for a smaller $r$ if necessary yet independent of $i$.
Note that this means that $F_{i+1}(D(r))$ is in the domain of
$\Phi_{i+1}$.  Also, we may choose this $r$ so that $F_{i+1}$ is a
biholomorphism onto its image.  By \textbf{(a3)} one can show that
these maps are essentially contracting in the sense that there are
$n_i\uparrow\infty$ such that
\begin{equation}\label{contracting1}
F_{n_{i+1}}\circ \cdot \cdot F_{n_i +2}\circ F_{n_i
+1}(D(r))\subset D(\frac {r}2)
\end{equation}
for every $i$. Now let us suppose each $F_i$ can be extended to a
biholomorphism of $\cn$, which we still denote as $F_i$, on $\cn$.
We may then let
$$
S_i=\lf(F_{n_i}\circ\cdots\circ F_2\ri)^{-1}(D(r))
$$
and $\Omega=\bigcup_iS_i$. Then $\Omega$ will be a pseudoconvex
domain in $\cn$ and it will be homeomorphic to $\mathbb{R}^{2n}$
because $S_i\subset S_{i+1}$. Let $\Psi$ be defined as
$$
\Psi(z)=\Phi_{n_i}\circ F_{n_i}\circ\cdots\circ F_2(z)
$$
for $z\in S_i$. Note that $$\Phi_{n_{i+1}}\circ
F_{n_{i+1}}\circ \cdots \circ F_{n_i}\circ\cdots \circ
F_2(z) =\Phi_{n_i}\circ F_{n_i}\circ\cdots\circ F_2(z)
$$
if $z\in S_i$ and thus $\Psi$ is a well-defined nondegenerate map from
$\Omega$ to $M$.  Moreover, using \textbf{(a3)} and that $g_i$ is `shrinking',
one can prove that $\Psi$ is surjective thus proving Theorem \ref{ChauTam-Stein2}.  In case $M$
is simply connected, one can also prove that $\Psi$ is
injective.

Now in general, $F_i$ can not be extended to a biholomorphism of
$\cn$.  The key is to show that the maps $F_i:D(r)\to \cn$ can be
approximated well enough by  biholomorphisms of $\cn$.  For this
one uses a theorem of Anderson-Lempert \cite{A} which states that
if $F$ is a biholomorphism from a star-shape domain in $\cn$ onto
a Runge domain in $\cn$, then $F$ can be uniformly approximated by
biholomorphisms of $\cn$ on compact subsets of the domain.  To use
the result it is thus sufficient to show that the image of
$F_i:D(r)\to \cn$ is Runge.  Once it is established that $F_i$ can
be approximated in this way, one uses the approximations to
construct a map $\Psi$ as above having the desired properties.

 From the above proof, one expects to obtain stronger results depending on
how much more can be said of the maps $F_i$.  We will see this to be the case
in the following sections \S5 and \S 6.

\section{Gradient \KR solitons}
In very special cases, the maps $F_i$ from the previous section will all
 be equal to a single map $F$.  This is the case when $g(t)$ is a gradient \KR soliton.

Let us recall that \KR solitons are solutions to the \KR
flow for which the metric evolves only by dilation and pull
back along a one parameter family of biholomorphisms.  More
specifically, $(M, g_{i \bar j})$ is said to be a \KRS if
there is a family of biholomorphisms $\phi_t$ on $M$, given
by a holomorphic vector field $V$, such that
$\G(x,t)=\phi_t^*(\G(x))$ is a solution of the normalized
\KR flow:

\begin{equation}\label{KRF1}
\begin{split}
\frac{\p}{\p t}&g_{i \bar j}(x, t)=-R_{i\bar j}(x, t)-\kappa g_{i \bar j}(x, t)\\
&g_{i \bar j}(x,0)=g_{i \bar j}(x)
\end{split}
\end{equation}
for some constant $\kappa$. The soliton is said to be {\it steady}
if $\kappa=0$, and {\it expanding} if $\kappa>0$ (which will be
normalized to be 1). In particular, note that $\phi_t:(M, g(t))
\to (M, g)$ is an isometry for every $t$, and thus \KR solitons
can be viewed as generalized fixed points for the \KR flow. If in
addition, the holomorphic vector field $V$ is the gradient of a
real valued function $f$, then the soliton is a gradient \KR
soliton and the metric $g_\ijb$ satisfies
\begin{equation}\label{gkrseqns}
\begin{split}
f_{i\bar j}&= R_{i\bar j}+\kappa  g_{i\bar j},\\
&f_{ij}=0.
\end{split}
\end{equation}
Conversely, if these equations are satisfied with bounded
curvature, then the corresponding solution to (\ref{s1e1}) is a
gradient \KR soliton by the uniqueness theorem of Chen-Zhu
\cite{Chen-Zhu2006}. These equations are thus known as gradient \KR
soliton equations.

 Solitons are extremely important in studying the
 formation of singularities under the flow and hence the underlying structure of the manifolds.
 As described in \cite{Ha4} by Hamilton, solitons
 typically arise as limit solutions when one takes a dilation limit
 around a singularity forming under the flow (see \cite{Ha2} on
 limit solutions to the Ricci flow).  In \cite{cao},
 Cao proved the following classification of limit solutions
 to the \KR flow.
\begin{thm}\label{caolimits}\cite{cao}
Let $(M, g(t))$ be a family of complete K\"ahler metrics on a
noncompact complex manifold $M$ which form  a solution to
(\ref{KRF1}) for $t\in[0,\infty)$ if $\kappa>0$ and for
$t\in(-\infty,\infty)$ if $\kappa=0$ such that the holomorphic
bisectional curvature is nonnegative and the  Ricci curvature is
positive. Assume that the scalar curvature of $g(t)$ assumes its
maximum in space time. Then $(M, g(t))$ is a gradient \KR soliton
which is steady if $\kappa=0$, and expanding if $\kappa=1$.
\end{thm}

The proof of Theorem \ref{caolimits} relies on the LYH type
Harnack inequality Theorem \ref{caoharnack} of Cao \cite{cao}.
The theorem suggests that studying the uniformization of gradient
\KR solitons should give insight into the uniformization of more
general manifolds as in \S6.

Let $\phi_t$ be the family of biholomorphisms defining a \KR
soliton so that $g(t)=\phi_t^*(g(0))$.  Assume the holomorphic
bisectional curvature is bounded and nonnegative. Suppose $p$ is a
fixed point of the flow $\phi_t$ and let $F=\phi_1=\phi_t|_{t=1}$.
Then the injectivity radii of $p$ with respect to $g(t)$ are
constant. Hence on can construct biholomorphism  $\Phi_i$ for
$g(i)$ as in the proof of Theorem \ref{ChauTam-Stein1}. In this
case, we can actually take $\Phi_i$ such that
$\Phi_{i+1}^{-1}=\Phi_1^{-1}\circ(\phi_1)^i$ and thus $
\Phi_{i+1}^{-1}\circ\Phi_i=\phi_1=F$. In other words, the $F_i$'s
from \S 4 can in this case be taken as the single map $F$.  And
studying the $F_i$'s reduces here to studying the discrete complex
dynamical system generated by $F$. Recall that $p$ is a fixed
point for $F$. In \cite{RR} Rosay-Rudin proved the following
result for {\it attractive} {\it basins} of a fixed points for
biholomorphisms of $\cn$ , which was pointed out to hold on
general complex $M$ by Varolin \cite{Varolin2000}:

\begin{thm}\label{RosayRudin}\cite{RR}
Let F be a biholomorphism from a complex manifold $M^n$ to itself
and let $p\in M^n$ be a fixed point for F.  Fix a complete
Riemannian metric g on M and define the basin of attraction
$$\Omega:=\{x\in M: \lim_{k\to \infty} dist_g (F^k(x),
p)=0\}$$where  $F^k=F\circ F^{k-1}, F^1=F$.

Then  $\Omega$ is biholomorphic to $\ce^n$ provided
$\Omega$ contains an open neighborhood around p.
\end{thm}

 Now in many cases, the dynamical system generated by $F=\phi_1$
 above can be shown to have a unique attractive fixed point with the whole manifold as a basin of attraction.  By Theorem
\ref{RosayRudin}, such a soliton must then be biholomorphic to
$\ce^n$.  This was observed by the authors  in \cite{CT2} where
they   proved:

\begin{thm}\label{t}\cite{CT2}
If $(M,g_{i\jbar})$ is a complete non-compact gradient
K\"ahler-Ricci soliton which is either steady with positive Ricci
curvature so that the scalar curvature attains its maximum at some
point, or expanding with non-negative Ricci curvature, then $M$ is
biholomorphic to $\ce ^n$.
\end{thm}

 This result was obtained independently by Bryant \cite{Bryant} for the steady
case.

To  prove the theorem one only needs to check that there is a unique fixed point of the biholomorphisms $\phi_t$ and that $\phi_t$ is
{\it contracting} on $M$.  These are easily verified using the
 condition on the positivity or nonnegativity of the Ricci
curvature.

  The proof of Theorem \ref{RosayRudin} relies on the fact that $F$
   can be transformed to have a normal form around $p$.
   This   fact is due to Sternberg \cite{St} for real systems, and was later independently proved by Rosay-Rudin for complex systems.  We sketch the proof here for the case $M= \ce^n$.  Let $F:\cn\to\cn$ be a biholomorphism such that
$F(0)=0$.  One then modifies $F$ by a biholomorphism
$T$ near the origin, so that (i) $T\circ F\circ T^{-1}$ is
close to an {\it upper triangular} map $G$; and (ii)
$T'(0)=I$, the identity map. Here $G=(g_1,\dots,g_n)$ is an
upper triangular map, if $g_1(z)=c_1z_1$,
$g_2(z)=c_2z_2+h(z_1)$,$\dots$,
$g_n(z)=c_nz_n+h(z_1,\dots,z_{n-1})$ for some constants
$c_1,\dots,c_n$.  By (i) we mean that for any $m$, we can choose $T$ and $G$ so that  $G^{-1}\circ T\circ F-T=O(|z|^m)$,
with $G$ being independent of $m$ when $m$ is large enough.
It can then by shown that $\Psi=\lim_{k\to\infty}G^{-k}\circ T\circ F^k$ exists
and is a biholomorphism from the basin onto $\cn$.  For any
$z$ in the basin, $F^k(z)$ is defined if $k$ is large
enough. On the other hand, $F$ is shrinking so $G$ is
expanding. Hence one may expect the image will be the whole
$\cn$.

In special cases, the \KR flow actually performs this
uniformization and the soliton metric converges in the re-scaled
subsequence sense to a complete flat \K metric under the flow. The
following result was obtained in  \cite{CT}.

\begin{thm}\label{soliton}\cite{CT} Let $(M, g_{i\jbar})$ be a complete non-compact
gradient \KRS as in Theorem \ref{t} with smooth potential
$f$ and equilibrium point $p$. let $g_{i\jbar}(x, t)$ be
the corresponding solution to (\ref{KRF1}) and let $\bold
v_p \in T^{1,0}_p (M)$ be a fixed nonzero vector with
$|\bold v_p|_0=1$. Then for any sequence of times
${t_k}\rightarrow \infty$, the sequence of complete \K
metrics $ \frac{1}{ |\bold v_p|_{t_k}^{2} } g_{i\jbar}(x,
t_k)$ subconverges on compact subsets of $M$ to a complete
flat \K metric $h_{i\jbar}$ on $M$ if and only if
$R_{i\jbar}(p)=\beta g_{i\jbar}(p)$ at $t=0$ for some
constant $\beta$.  In particular,if this condition is
satisfied then $M$ is biholomorphic to $\ce^n$.
\end{thm}

 Note that the Theorem suggests in general, we do not expect to prove uniformization by rescaling a solution $g(t)$ to the \KRF to obtain a complete K\"ahler flat metric as a limit.

 \section{Eternal solutions to the normalized \KRF}
 In this section we generalize Theorem \ref{t} for \KR solitons to
 $eternal$ solutions to the normalized \KR flow.  As a corollary of this
 we will present a uniformization theorem for the case of average quadratic curvature decay.
This is a critical case for uniformization in light of the gap
phenomenon for manifolds with faster than quadratic curvature
decay discussed in \S 7.  The first major result in this case was
the following theorem of Mok \cite{M1}.
\begin{thm}\label{Mok}\cite{M1}
Let $(M^n, g)$ be a complete noncompact \K surface with positive
holomorphic bisectional curvature.  Suppose that the following
conditions are satisfied for some $p \in M^n$,
\begin{enumerate}
\item[(i)] $Vol(B_p(r))\geq C_1 r^{2n}$, for all $r\geq 0$
\item[(ii)] $\mathcal{R}(x)\leq \frac{C_2}{(d(p, x)
+1)^{2}}$. for some $C_1, C_2, \e>0$.
\end{enumerate}
Then $M^n$ is biholomorphic to an affine algebraic variety.
\end{thm}
It is well known that a complete noncompact Riemannian manifold
with positive sectional curvature is diffeomorphic to the
Euclidean space \cite{GromollMeyer69}. Using a result of Ramanujam
\cite{Ramanujam1971} which states that an algebraic surface which
is homeomorphic to $\mathbb{R}^4$ must be biholomorphic to $\C^2$,
Mok \cite{M1} concluded:
\begin{cor}\label{Mok-cor}
Let $(M^2,g)$ be a complete noncompact \K surface with positive
sectional curvature satisfying conditions (i) and (ii) in the
above theorem. Then   $M$ is biholomorphic to $\C^2$.
\end{cor}

 The method of Mok is to
 construct enough polynomial growth holomorphic functions to embed $M$ into
 some $\C^N$ so that the image will be an affine algebraic variety.
 The proof uses algebraic geometric methods.
   Later, in separate works, Chen-Tang-Zhu, Chen-Zhu and Ni used the
\KRF to improve Mok's result in Theorem  \ref{Mok}. In particular
 in the case of $n=2$,  Chen-Zhu \cite{CTZ,
Chen-Zhu05} obtained the same result as in the Corollary
\ref{Mok-cor} by only assuming   positive and bounded bisectional
curvature and maximal volume growth. Ni \cite{Ni05} further
weakened the condition of positive holomorphic bisectional
curvature to nonnegative holomorphic bisectional curvature.
 The main idea is to show that maximum volume
growth still implies quadratic curvature decay condition, as
mentioned in Theorem \ref{Ni-curvature}.  Then one can still prove
that $M$ is homeomorphic to $\mathbb{R}^4$ and produce enough
polynomial growth holomorphic functions to carry over Mok's method
and result of Ramanujam. Chen-Zhu \cite{Chen-Zhu02} also proved
that if a K\"ahler surface with bounded and positive sectional
curvature is such that the integral of $(Ric)^2$ is finite,   then
the surface is biholomorphic to $\mathbb{C}^2$, using Ramanujam's
result again.  Ramanujam's theorem however is only true for
complex surfaces, and for higher dimensions we need other methods.

 We would like to use the results on \KR flow by Shi in \S2 to
  generalize Theorem \ref{t} for \KR solitons to general
  solutions $g(t)$ to the \KR flow.  It is natural here
  to consider $eternal$ solutions to (\ref{KRF1}),
  in other words solutions defined for $t\in (-\infty, \infty)$.
  It is readily seen that a steady or expanding gradient
  \KR soliton is indeed an eternal solution to (\ref{KRF1}).
  In light of this, one may expect that Theorem \ref{t} is still
  true when $g(t)$ is an eternal solution to (\ref{KRF1})
  with nonnegative uniformly bounded holomorphic bisectional curvature.
  This expectation is confirmed in Theorem \ref{CTeternal}
  and was proved by the authors in \cite{CT4}.
  Before stating the theorem, we first discuss the case of
  quadratic curvature decay in relation to eternal solutions.

 Suppose $(M,g)$ is complete noncompact with
bounded and nonnegative holomorphic bisectional curvature so that
its scalar curvature satisfies the quadratic decay condition
(\ref{quadraticdecay1}).  Then by the results in \S 2,
(\ref{KRF2})  has a long time solution $\tg$ with initial data
$g$, and the scalar curvature $\mathcal{\wt R}$ will satisfy $t\wt
R\le C$ for some constant $C$ uniform on spacetime.  Now if we let
$g(t)=e^{-t}\tg(e^t)$, then $g(t)$ will be an eternal solution to
the normalized \KR flow (\ref{KRF1}) for $\kappa=1$.  Moreover, it
is easy to see that $g(t)$ has uniformly bounded nonnegative
holomorphic biesctional curvature.  Conversely, given such an
eternal solution $g(t)$ one sees that $\ttg(t)=tg(\log t)$ solves
the unnormalized \KR flow
\begin{equation}\label{KRF2}
    \frac{\p}{\p t}\ttg_{\ijb}=-\wt R_{\ijb}.
\end{equation}
for $t\ge1$, and that $t\wt R\le C$ for some uniform constant $C$.
This in turns implies that the scalar curvature of the initial
metric  $g(0)$ satisfies the quadratic decay condition
(\ref{quadraticdecay1})  see \cite{NT,Ni05}.

Now we state our theorem on eternal solution:
\begin{thm}\label{CTeternal}\cite{CT4} Let $M^n$ be a noncompact complex manifold.
Suppose there is a smooth family of complete \K metrics $g(t)$ on
$M$ such that for $\kappa=0$ or $1$, $g(t)$ satisfies
\begin{equation}\label{s1e3}
\frac{\p }{\p t} {g}_{i\jbar}(x,t)=-{R}_{i\jbar}(x,t)-\kappa
g_\ijb(x,t)
\end{equation}
for all $t\in (-\infty,\infty)$ such that for every $t$, $g(t)$
has uniformly bounded non-negative holomorphic bisectional
curvature on $M$  independent of $t$. Then $M$ is holomorphically
covered by $\ce^n$.
\end{thm}
 By the remarks preceding the theorem, we have the following
  result by the authors:

\begin{thm}\label{CTQuadDecay}\cite{CT4}
Suppose $(M^n, g)$ has holomorphic bisectional curvature
which is bounded, non-negative and has average quadratic
curvature decay.
 Then $M$ is holomorphically covered by
$\ce^n$.
\end{thm}
Combining this with Theorem \ref{Ni-curvature}, we conclude:

\begin{cor}\label{CTMaxVol} Let $(M^n,g)$ be a
complete noncompact K\"ahler manifold with bounded and nonnegative
holomorphic bisectional curvature such that $M$ has maximum volume
growth  then $M$ is biholomorphic to $\cn$.
\end{cor}

\begin{rem}
Corollary \ref{CTMaxVol} was proved before Theorem
\ref{CTQuadDecay} by the authors in \cite{CT3}.
\end{rem}

As noted earlier, if we assume the holomorphic bisectional
curvature is bounded and nonnegative, then Theorem \ref{CTeternal}
is basically a direct generalization of Theorem \ref{t} for
gradient \KR solitons. However, the proof of Theorem
\ref{CTeternal} is much more complicated.  Beginning with a
solution $g(t)$ to the \KRF as in Theorem \ref{CTeternal}, fix
some point $p\in M$ and construct maps $\Phi_i$ as in the proof of
Theorem \ref{ChauTam-Stein2} in \S 4.  For simplicity, we will
assume that $\Phi_t$ is injective for every $t$ (thus $M$ is
simply connected).  In other words, we assume the injectivity
radius of $g(t)$ at $p$ is bounded from below independently of
$t$.  Such a bound exists in the case of \cite{CT3}, where maximum
volume growth is assumed and removing the dependence on this bound
is the essential generalization made in \cite{CT4}. We may also
assume that $g(t)$ has positive Ricci curvature because of a
dimension reduction result of Cao \cite{cao1}.

 Now
for $N>0$ sufficiently large, as in \S4, we can find a
sequences of biholomorphisms
 $F_i$  from $D(r)$ onto its image which is inside $D(r)$:
\begin{equation}\label{Fi's}
F_{i+1}=\Phi_{(i+1)N}^{-1}\circ \Phi_{iN}:D(r)\to D(r)
\end{equation}
for $i\geq 1$. These $F_i$'s are basically the same as those in
the proof of Theorem \ref{ChauTam-Stein2}, which as  noted in \S
5, can be chosen to be a single map $F$ when $g(t)$ is gradient
\KR soliton. One would now like to imitate the proof of
Rosay-Rudin's Theorem \ref{RosayRudin}.  A key step in their proof
was to transform $F$ into a particularly nice form.  Now the main
difficulty here is simultaneously transforming the sequence
$\{F_i\}$ into a likewise nice form. In \cite{JV},
Johnsson-Varolin showed that this can be done provided that asymptotically they behave close enough
to a single map $F$.  This closeness is essentially in terms of
the Lyapunov regularity of the $F_i$'s  (see \cite{BP} for the
terminology).  In terms of the \KR flow, the authors proved
\cite{CT4} that this transformability is possible due to the
Lyapunov regularity of $g(t)$ as described in the following:

\begin{thm}\label{Lyapunov}\cite{CT4} Let $M^n, g(t)$ be as in Theorem \ref{CTeternal}
such that the Ricci curvature of $g(t)$ is positive.  Let $p\in M$
be fixed and let $\lambda_1(t)\ge \dots\ge\lambda_n(t)>0$ be the
eigenvalues of $R_\ijb(p,t)$ relative to $g_\ijb(p,t)$. Then
$\lim_{t\to\infty}\lambda_i(t)$ exists for $1\le i\le n$ and
  there is a constant $C>0$ such that $\lambda_n(t)\ge C$ for
all $t\ge0$. Moreover, if we let $\mu_1>\dots>\mu_l> 0$ be the
distinct limits   with multiplicities $d_1,\dots,d_l$, then
$V=T_p^{(1,0)}(M)$ can be decomposed orthogonally with respect to
$g(0)$ as
 $V_1\oplus \cdot\cdot\cdot\oplus V_l$ so that the following are true:

\begin{enumerate}
 \item [(i)] If $v$  is a nonzero vector in  $V_i$ for some $1\le
i\le l$, then
 $\lim_{t \to \infty} |v_i(t)|=1$ and thus  $\lim_{t\to\infty}Rc(v(t),\bar v(t))=\mu_i$ and
$$\lim_{t \to \infty} \frac{1}{t}\log \frac{|v|_t^2}{|v|_0^2}=-\mu_i-1.$$  Moreover, the convergences are uniform over all
 $v\in V_i\setminus\{0\}$.
 \item[(ii)] For $1\le i, j\le l$ and for nonzero vectors  $v \in
V_i$ and $w \in V_j$ where $i\neq j$, $\lim_{t\to
\infty}\langle v(t),w(t)\rangle_t=0$ and the convergence is
uniform over all such nonzero vectors $v, w$. \item [(iii)]
$\dim_\C(V_i)=d_i$  for each $i$.
 \item[(iv)]
$$\sum_{i=1}^l(-\mu_i-1)\dim_\C V_i=\lim_{t\to\infty}
\frac1t\log \frac{\det(g_{i\bar j}(t))}{\det ({g}_{i\bar
j}(0)}.$$
\end{enumerate}
\end{thm}
The theorem basically says that the eigenvalues of the Ricci
tensor have limits and the eigenspaces are almost the same.

The proof of the theorem relies on an important differential
Li-Yau- Hamilton(LYH) inequality for the \KRF by Cao  \cite{cao}:

\begin{thm}\label{caoharnack}\cite{cao}
Let $g(t)$ be a complete solution to the \KR flow (\ref{s1e1}) on
$M\times [0, T)$ with bounded and nonnegative holomorphic
bisectional curvature. For any $p\in M$ and holomorphic vector $V$
at $p$, let
$$
Z_{i\bar{j}}=\frac{\p}{\p t}
R_{i\bar{j}}+R_{i\bar{k}}Z_{k\bar{j}}+R_{i\bar{j},
k}V_{\bar{k}}+R_{i\bar{j}, \bar{k}}V_{k}
+R_{i\bar{j}k\bar{l}}V_{\bar{k}}V_l+\frac{1}{t}R_{i\bar{j}}.
$$
Then $$Z_{i\bar{j}}W^i W^{\bar{j}}\ge0$$ for all
holomorphic vectors $W$ at $p$.
\end{thm}
Using this differential inequality the authors proved that
\begin{thm}\label{TH}\cite{CT5}
Let $g(t)$ be a complete solution to (\ref{s1e1}) with
non-negative holomorphic bisectional curvature such that for any
$T>0$, $g(t)$ has bounded curvature for all $t\in[0,T]$. Fix some
$p \in M$ and let $\lambda_i(t)$ be the eigenvalues of $Rc(p, t)$
arranged in increasing order. Then
$$t\lambda_k(t)$$ is nondecreasing in  $t$ for all $1\leq k \leq n$.
\end{thm}

Now under the condition of Theorem \ref{Lyapunov}, the proof of Theorem \ref{TH} directly implies that
$\lambda_i(t)$ is nondecreasing in $t$ for every $1\leq i \leq n$. This will imply that
$\lim_{t\to\infty}\lambda_i(t)$ exists for all $i$.  From
this, one argues as in the proof of Theorem \ref{caolimits}
in \cite{Cao} to prove that $g(t)$ behaves like gradient
\KR soliton with fixed point at $p$ as $t\to\infty$ in the
following sense: For any $t_k\to\infty$, there is a
subsequence of $g(t+t_k)$ such that $(M,g(t+t_k))$ converge
to a gradient \KR soliton.  To prove this one actually only needs
the convergence of the scalar curvature $R(p,t)$. In case
the manifold has maximum volume growth, a more general
result similar to this was obtained by Ni \cite{Ni005}
independently. Now it is easy to see that if $g(t)$ is a
\KR soliton with fixed point $p$ then Theorem
\ref{Lyapunov} is true, and that in this case we do not
even have to take limits.  Observing this, one then argues carefully to
obtain the results in Theorem \ref{Lyapunov}.

 We now return to our sketch of proof of Theorem \ref{CTeternal}.
 The $F_i's$  define a randomly iterated
complex dynamical system on $D(r)$ with fixed point at the
origin.  Moreover, using that the
 Ricci curvature is bounded away from zero at $p$ for all $t$ by Theorem \ref{Lyapunov}, one can show that the maps $F_i$ are uniformly
contracting at the origin and that $\lim_{i\to \infty} F_i\circ
\cdot \cdot\circ F_1 (D(r))=0$.  This is one of the biggest
differences between the $F_i$'s here and those in \S
4.  Although the maps there were eventually contracting, they are
in general not uniformly contracting.

 Theorem \ref{Lyapunov} can now be translated into Lyapunov regularity of the
system $\{F_i\}$ which can roughly be described as follows.  Let $A_i=F_i'(0)$.  Then the system is Lyapunov regular at $0$ if one can decomposed $\cn$ orthogonally with respect the Euclidean metric as $E_1,\dots,E_l$ such that if $E_k^{(i+1)}=A_{i+1}(E_k^{(i)})$, then $E_k^{(i)}$ are asymptotically
orthogonal and for each $k$,  $A_i$ is asymptotically contracting at a constant rate on $E_k^{(i)}$.

Once this is established, one follows the construction in \cite{JV}
to construct a sequence of biholomorphisms $G_i:\ce^n \to
\ce^n$, and $T_i:D(r) \to D(r)$ (it  may be
necessary to take $r$ smaller here, but independently of
$i$).   Here, the $G_i$'s represent approximations of the
$F_i$'s in $\text{Aut}(\ce^n)$ which are lower triangular
in a certain sense while the $T_i$'s represent a sequence
of holomorphic coordinate changes of $D(r)$. The following
Lemma describes the extent of this
approximation (\cite{CT3} Lemma 5.7).

\begin{lem}\label{corrections} Let $k\ge0$ be an integer. Then
$$
\Psi_k=\lim_{l\to\infty}G^{-1}_{k+1}\circ G^{-1}_{k+2}\circ\cdots\circ
G^{-1}_{k+l}\circ T_{k+l}\circ F_{k+l}\circ\cdots\circ F_{k+2}\circ F_{k+1}
$$
exists and is a nondegenerate holomorphic map from $D(r)$
into $\C^n$. Moreover, there is a constant $\gamma>0$ which
is independent of $k$ such that
$$
\gamma^{-1}D(r)\subset\Psi_k (D(r))\subset \gamma D(r).
$$
\end{lem}

Since $F_i$ is uniformly contracting, one can show that  that the
sets $\Omega_i=\Phi_{iT}(D(r))$ exhaust $M$ as $i\to \infty$. This
together with Lemma \ref{corrections} and the definition of the
$F_i's$ tell us that the sequence of maps
$$
S_i=G^{-1}_1\circ\cdots G^{-1}_{i}\circ T_{i}\circ \Phi_{iN}^{-1}:
\Omega_i \to \cn
$$
converges to a biholomorphic map $\Psi$ from $M$ into
$\ce^n$.   It is shown in (\cite{CT3}, \S5) that $\Psi$ is
onto, and thus $M$ is biholomorphic to $\ce^n$.

Now in case that the injectivity radius of $p$ with respect to $g(t)$
is not bounded away from zero, one works on the pullback
metrics $\hat g (t)$ of $g(t)$ under the exponential maps.  In this setting the
injectivity radius will be bounded away from zero and one can
show that $\hat g (t)$ still behaves like gradient \KR soliton locally
near $0$ as $t\to\infty$.  One then constructs maps $F_i$ as in
the proof of Theorem \ref{ChauTam-Stein1}, which one then shows to be Lyapunov regular, and proceeds as above to obtain Theorem \ref{CTeternal}.

The fact that $\hat g$ behave like gradient \KR solitons locally near $0$ as
$t\to\infty$ can be used to prove the following corollary to the
Theorem \ref{CTeternal}.

\begin{cor}\label{CTeternal-cor} Let $M, g(t)$ as in Theorem
\ref{CTeternal}. Suppose the Ricci curvature is positive with
respect to $g(0)$ at some point $p$. Then $M$ is simply connected
and is biholomorphic to $\cn$. In particular, if $(M,g)$ is a
complete noncompact \K manifold with bounded positive holomorphic
bisectional curvature which satisfies the quadratic decay
condition (\ref{quadraticdecay1}), then $M$ is biholomorphic to
$\cn$.
\end{cor}
One might want to compare the last assertion of the corollary with
the statement of Yau's conjecture and the result of Zheng's
Theorem \ref{Zheng}.  The main point of the proof of the corollary
is to show that the first fundamental group of $M$ is finite.  Since $M$ is covered by $\cn$, $M$ must then
be simply connected by a well-known result, see \cite{Brown1982}.

 \section{A Theorem of Mok-Siu-Yau and its generalizations}
One may expect stronger results when the curvature decays faster than quadratic.  In fact, there are  gap theorems which
tell us that curvature of a nonflat complete noncompact K\"ahler
manifold $(M, g)$ with nonnegative holomorphic bisectional
curvature cannot decay too fast.  These can be viewed as
converses to the curvature decay Theorem \ref{ChenZhu-curvature}.
The following classic gap theorem of Mok-Siu-Yau \cite{MSY} in
1981 was the first result supporting Yau's conjecture.

\begin{thm}\label{MSY-gap}\cite{MSY}
Let $M$ be a complete noncompact \K manifold with nonnegative
holomorphic bisectional curvature.  Suppose that the following
conditions are satisfied for some $p \in M$ and some $\epsilon>0$:
\begin{enumerate}
\item[(i)] $Vol(B_p(r))\geq C_1 r^{2n}$, for all $r\geq 0$.
\item[(ii)] $\mathcal{R}(x)\leq \frac{C_2}{(d(p, x)
+1)^{2+\epsilon}}$. for some $C_1, C_2, \e
>0$. \item[(iii)] Either $M$ is Stein or $M$ has nonnegative
sectional curvature.
\end{enumerate}
Then $M$ is isometrically biholomorphic to $\C^n$.
\end{thm}
By Theorem \ref{NiTam-Stein}, condition (iii) is superfluous
because of (i).  The condition of maximum volume growth (i) however
seems rather strong. In   \cite{Chen},  Chen-Zhu proved the
following:

 \begin{thm}\label{Chen-Zhu-gap}\cite{Chen}
Let $M$ be a complete noncompact \K manifold with nonnegative and
bounded holomorphic bisectional curvature.  Suppose that for some
positive real function $\epsilon(r)$ satisfying $\lim_{r\to
\infty}\epsilon (r)= 0$ we have
\begin{equation}\label{faster-quadratic}
\frac{1}{Vol(B_{x}(r))}\int_{B_{x}(r)} \mathcal{R} dV\leq
\frac{\epsilon(r)}{(1+r)^2}
\end{equation}   for all  $x\in M$ and for all $r>0$.   Then the
universal cover of $M$    is isometrically biholomorphic to
$\C^n$.
\end{thm}
The theorem says that if $M$ has nonnegative holomorphic
bisectional curvature such that   the curvature is bounded
and decays faster than quadratic on average   {\it
uniformly} at all points in $M$, then $M$ is flat. It is
not hard to see these hypothesis are implied by (i) and (ii)
in Theorem \ref{MSY-gap}.  The proof of Theorem
\ref{Chen-Zhu-gap} uses the \KR flow.  In fact, one can see
from  the proof of Theorem \ref{Shi-longtime} that if we
consider the normalized flow (\ref{KRF1}) with $\kappa=1$,
then the scalar curvature will tend to zero at infinity.
Hence it must be zero initially by the monotonicity derived
from Theorem \ref{caoharnack}. Note that in order to use
the \KR flow, one needs to assume  that the curvature is
bounded.  On the other hand, there is a way to prove a
stronger result in this case without using the \KR flow by
modifying the method of \cite{MSY}.

Let us go back to the proof of Mok-Siu-Yau \cite{MSY}.
Their method  is as follows: First one solves the Poisson
equation $\frac12\Delta u=\mathcal{R}$ with good estimates.
This can be done because  assumptions (i) and (ii) in
Theorem \ref{MSY-gap} give a good estimate of the Green's
function. Secondly, one can show that $||\ii \p\bar\p
u-\Ric||$ is subharmonic using the fact that $M$ has
nonnegative holomorphic bisectional curvature. Using the
estimate of $u$, one may get an integral estimate of
$||\p\bar\p u||^2$ on geodesic balls. Then by a mean value
inequality, we conclude that $\ii \p\bar\p u=\Ric$. In
particular, $u$ is plurisubharmonic.  Finally, one proves
that $u$ is constant implying that $M$ is flat.  One proves
this by contradiction: assuming $u$ is not constant, one
produces a function $v$ which at a point is strictly
plurisubharmonic and satisfies $(\ii \p\bar\p v)^n\equiv0$.
If the manifold is Stein one embeds $M$ in $\C^N$ for some
$N$ and proceeds to use the coordinate functions of $\C^N$
to construct such a $v$.  If $M$ has nonnegative sectional
curvature one uses the Busemann function, together with
$u$, to construct $v$.

Following this line of argument, Ni-Shi-Tam in \cite{NST} obtained
a general result on the existence of solution of the Poisson
equation on complete noncompact Riemannian manifold with
nonnegative Ricci curvature without any volume growth condition.
   In particular, it was shown there that if
(\ref{faster-quadratic}) holds at some point $p$, then one can still
solve:

 $$
 \frac12\Delta u=\mathcal{R}
 $$
 with
 \begin{equation}\label{sub-log}
    \limsup_{r\to\infty}\frac{u(x)}{\log d(p,x)}\le0.
\end{equation}
These follow from the classic results of Li-Yau
\cite{LiYau86} on the heat kernel estimates on manifolds
with nonnegative Ricci curvature, which provide an
estimates for the Green's function without assume maximum
volume growth.

In order to prove that $\ii \p\bar\p u=\Ric$, one needs the
condition on
 the $L^2$ norm of $\mathcal{R}$:
 \begin{equation}\label{L2-norm}
   \liminf_{r\to\infty}\lf[\exp(-ar^2)\int_{B_p(r)}\mathcal{R}^2
   dV\ri]<\infty
\end{equation}
for some $a>0$. Note that  $\mathcal{R}$ may be allowed to
grow like $\exp(a'r^2)$ for some $a'>0$.  From this, one
may then solve the heat equation with initial data $u$ and
use the maximum principle and the mean value inequality as
developed by Li-Schoen \cite{LiSchoen84} to conclude that
$u$ is indeed a potential of the Ric form. Here
(\ref{L2-norm}) is used to apply the maximum principle.

Finally, one can prove that $u$ is constant by (\ref{sub-log}) and
the following Liouville property for pluirsubharmonic functions of
Ni-Tam \cite{NT2}:
\begin{thm}\label{Ni-Tam-Liouville}\cite{NT2}
Let $(M,g)$ be a complete noncompact K\"ahler manifold with
nonnegative holomorphic bisectional curvature. Suppose $u$
is a continuous plurisubharmonic function satisfying
(\ref{sub-log}), then $u$ is constant.
\end{thm}
We may assume here that $u$ is bounded from below. The idea of proof of Theorem \ref{Ni-Tam-Liouville} uses methods as in the proof of Theorem
\ref{NiTam-curvature} and uses the following result of Ni
\cite{Ni} to conclude that $u$ is actually harmonic:

\begin{thm}\label{Ni-Liouville}\cite{Ni}
Let $(M^n,g)$ be a complete noncompact K\"ahler manifold
with nonnegative Ricci curvature. Suppose $u$ is a
plurisubharmonic function on $M$ satisfying
(\ref{sub-log}). Then $(\p\bar \p u)^n\equiv0$.
\end{thm}
Hence $u$ must be harmonic and therefor constant by a classical
result of Cheng-Yau \cite{ChengYau75} on harmonic function
on complete manifold with nonnegative Ricci curvature. From these results
Ni-Tam \cite{NT2} proved the following, which  is the best
result up to now in the generalization of Theorem
\ref{MSY-gap}:

\begin{thm}\label{Ni-Tam-gap}\cite{NT2}
Let $(M,g)$ be a complete noncompact \K manifold with
nonnegative holomorphic bisectional curvature.  Suppose the
scalar curvature $\mathcal{R}$ satisfies
(\ref{faster-quadratic}) and (\ref{L2-norm}) for some $p$.
Then the universal cover of $M$ is isometrically
biholomorphic to $\C^n$.
\end{thm}

It is still unknown whether the condition (\ref{L2-norm}) can be
removed.

 \bibliographystyle{amsplain}

\end{document}